\documentclass[12pt]{article}
\usepackage{fullpage}
\usepackage{amssymb, amsmath, amsfonts}
\newtheorem{theorem}{Theorem}[section]

\newtheorem{definition}[theorem]{Definition}
\newtheorem{proposition}[theorem]{Proposition}
\newtheorem{corollary}[theorem]{Corollary}

\title{Bounds for avalanche critical values of the Bak-Sneppen model}
\author{Alexis Gillett, Ronald Meester and Misja Nuyens\\ Vrije Universiteit 
Amsterdam\\ 
The Netherlands}
\date{}

\begin{document}

\maketitle

\begin{abstract}

We study the Bak-Sneppen model on locally finite transitive graphs $G$, in 
particular on  $\mathbb{Z}^d$ and on $T_{\Delta}$, the regular tree with common 
degree $\Delta$.  
We show that the avalanches of the Bak-Sneppen model dominate independent site 
percolation, in a sense to be made precise.  Since avalanches of the Bak-Sneppen 
model are dominated by a simple branching process, this yields upper and lower 
bounds for the so-called avalanche  critical value $p_c^{BS}(G)$. Our main results 
imply that  $\frac{1}{\Delta+1} \le p_c^{BS}(T_\Delta) \le \frac{1}{\Delta -1}$, 
and that $\frac{1}{2d+1}\leq p_c^{BS}(\mathbb{Z}^d)\leq \frac{1}{2d}+  
\frac{1}{(2d)^2}+O\big(d^{-3}\big)$, as $d\to\infty$.\\ 

\end{abstract}

\noindent
\section{Introduction and main results}

The Bak-Sneppen model was originally introduced as a simple model of evolution  by 
Per Bak and Kim Sneppen \cite{baksneppen}. Their original model can be defined as 
follows. There are $N$ species (vertices) arranged on a circle, each of which has 
been  assigned a random \textit{fitness}. The fitnesses are independent and 
uniformly  distributed on $(0,1)$. At each discrete time step the system evolves 
by  locating the lowest fitness and replacing this fitness, and those of its two  
neighbours, by independent and uniform $(0,1)$ random variables. We say that a 
vertex whose fitness is changed by this procedure has been {\em updated}.

It is not particularly significant that the underlying graph of the model  is the 
circle, or $\mathbb{Z}$ in the thermodynamic limit. Bak-Sneppen models  can be 
defined on a wide range of graphs using the same update rule as above: the vertex 
with minimal fitness  and its neighbours are updated.   Unlike particle systems 
such as percolation or the contact process, the Bak-Sneppen model has no tuning 
parameter. Therefore, it has been described as exhibiting self-organised critical 
behaviour, see  \cite{corrie} for a discussion.

One of the ways to analyse Bak-Sneppen models is to break them down into a series 
of {\em avalanches}. An avalanche from a threshold $p$, referred to as a 
$p$-avalanche, is said to occur between times $s$ and $s+t$ if at time $s$ all the 
fitnesses are equal to or greater than $p$ with at most one vertex where equality 
holds, and time $s+t$ is the first time after $s$ at which  all fitnesses are 
larger than $p$. The vertex 
with minimal fitness at time $s$ is called the  {\em origin} of the avalanche. A 
$p$-avalanche can be considered as a stochastic process in its own right. The key 
feature of the origin is that it has the minimal fitness (as it will be updated 
immediately). Hence, we can consider its fitness to be any value, as long as this 
value is minimal.  
Vertices with fitness below the threshold are called {\em active}, others are 
called {\em inactive}. Note that the exact fitness value of an inactive vertex is 
irrelevant for the avalanche, since this value can never be minimal during the 
avalanche. This motivates the following formal definition of an avalanche.

\begin{definition} A $p$-avalanche with origin $v$  on a graph $G$ (with vertex set 
$V(G)$)  is a stochastic process with  state space $\{[0,p]^A,A \subset V(G)\}$ and 
initial state $p^{\{v\}}$. The process follows the update rules of the Bak-Sneppen 
model. Any vertex with a fitness smaller than or equal to $p$ is included. Any 
vertex with a fitness larger than $p$ is not included. The process terminates when 
it is the empty set.
\end{definition}

Studying avalanches has considerable advantages. A Bak-Sneppen model on an {\em 
infinite} graph is not well-defined: when there are infinitely many  vertices, 
there may not be a vertex with minimal fitness. However, Bak-Sneppen {\em 
avalanches} can be defined on any locally finite graph as follows:  at time 0 all 
vertices have fitness 1, apart from one vertex, the  origin of the avalanche, which 
has fitness $p$. We then apply the update rules of the Bak-Sneppen model, until all 
fitnesses are above $p$.  This is consistent with our previous notion, as it is 
only the fitnesses updated during the avalanche that determine the avalanche's 
behaviour. The ability to look directly at infinite graphs is very desirable, 
because the most interesting  behaviour of the Bak-Sneppen model is observed in the 
limit as the number of vertices in the graph tends to infinity.

In the literature alternative types of avalanches have been proposed, see 
\cite{lee, li}.  The definition given here corresponds to  the most commonly used 
notion of an avalanche and was introduced by Bak and Sneppen \cite{baksneppen}.  
For a more thorough  coverage readers are directed to Meester and Znamenski 
\cite{zna1, zna2}. Note that  unlike the Bak-Sneppen model itself, the avalanches 
do have a tuning parameter, namely the threshold $p$.

In this paper, we look mainly at {\em transitive} graphs. The behaviour of an  
avalanche on a transitive graph is independent of its origin: an avalanche  with 
origin at vertex $v$ behaves the same as an avalanche  with origin 0. When 
analysing avalanches on transitive graphs, it is therefore natural to talk about a 
typical $p$-avalanche  without specifying its origin. 
To analyse avalanches, some definitions are needed. The set 
of vertices updated by  an avalanche is referred to as its range set, with the {\em 
range} being the  cardinality of this range set. Letting $r^{BS}_G(p)$ denote the 
range of a $p$-avalanche on a transitive graph $G$, we define the  {\em (avalanche) 
critical value} of the Bak-Sneppen model as
 \begin{equation}\label{critval}p_c^{BS}(G) = 
\inf\{p:\mathbb{P}(r^{BS}_G(p)=\infty)>0\}.\end{equation} 

Numerical simulations  \cite{baksneppen} suggest that the stationary marginal 
fitness distributions  for the Bak-Sneppen model on $N$ sites tend to a uniform 
distribution on   $(p^{BS}_c(\mathbb{Z}),1)$, as $N \to \infty$. It has  been 
proved in \cite{zna2} that this is indeed the case if 
$p_c^{BS}(\mathbb{Z})=\widehat{p}_c^{BS}(\mathbb{Z})$, where 
$\widehat{p}_c^{BS}(\mathbb{Z})$ is another critical value, based on the expected 
range, and is defined as 
\begin{equation}\label{crit2} 
\widehat{p}_c^{BS}(G)=\inf\{p:\mathbb{E}[r^{BS}_G(p)]=\infty\}.\end{equation}
It is widely believed, but unproven, that these two  critical values are equal.
 
It should now be clear that knowledge about the value of $p_c^{BS}(G)$ is vital in 
determining the self-organised limiting behaviour of the Bak-Sneppen model, even 
though there is no tuning parameter in the model.  Although in this paper we focus 
on the  critical value (\ref{critval}), our bounds for the critical value 
(\ref{critval}) also hold for the critical value (\ref{crit2}), see Section 
\ref{sec6}.

The approach of this paper is to compare Bak-Sneppen avalanches with two  
well-studied processes, namely branching processes and independent site 
percolation.  A simple comparison with branching processes gives a lower bound on 
the critical  value, whereas a  more complex comparison with site percolation gives 
an  upper bound. To warm up, we first give the (easy) lower bound. 

\begin{proposition}\label{thm1}  On any locally finite transitive 
graph $G$ with common vertex degree $\Delta$, we have
\[  p_c^{BS}(G) \ge \frac{1}{\Delta+1}.\]
\end{proposition}
\noindent
{\bf Proof:}  At every discrete time step of the system, we draw $\Delta + 1$ 
independent  uniform  $(0,1)$ random variables  to get the new fitnesses of the 
vertex with minimal fitness, and of its $\Delta$ neighbours. Each of these $\Delta 
+1$  new fitnesses is below the threshold $p$ with probability $p$,  independent of 
each other. This induces a coupling with a simple branching  process with binomial 
$(\Delta+1,p)$ offspring distribution, where every active  vertex in the 
Bak-Sneppen avalanche is  represented by at least one particle  in the branching 
process. Hence, if the branching process dies out, then so does  the Bak-Sneppen 
avalanche. Therefore the critical value of the Bak-Sneppen  avalanche can be no 
smaller than the critical value of the branching  process.~\hfill $\Box$\\

The main result of this paper is the following upper bound for the critical value  
$p_c^{BS}(G)$ of the Bak-Sneppen model on a locally finite transitive graph $G$.  
The critical value for independent site percolation on $G$ is denoted by 
$p_c^{site}(G)$.  We recall that for site percolation on $G$ with parameter $p$,  
the probability of an infinite cluster at the origin is positive for all  
$p>p_c^{site}(G)$, and 0 for all $p<p_c^{site}(G)$.

\begin{theorem}\label{thm2} On any locally finite transitive 
graph $G$, we have
\[  p_c^{BS}(G) \le p_c^{site}(G). \]
\end{theorem} 
This result implies that on many locally finite transitive graphs, $p_c^{BS}$ is 
non-trivial. For the Bak-Sneppen avalanche on $\mathbb{Z}$, Theorem \ref{thm2} 
gives a trivial upper bound, but in this case we know  from  \cite{zna1} that  
$p_c^{BS}(\mathbb{Z})\leq 1-\exp(-68)$ .

Since the critical value of site percolation on $T_{\Delta}$, the regular tree with 
common degree $\Delta$, equals $1/(\Delta-1)$, the following corollary holds.
\begin{corollary}
The critical value of the Bak-Sneppen model on a regular tree $T_\Delta$, with 
common degree $\Delta$, satisfies
\[ \frac{1}{\Delta+1} \le p_c^{BS}(T_\Delta) \le \frac{1}{\Delta -1}. \] 
\end{corollary}

Applying the expansion for the critical value of site percolation on $\mathbb{Z}^d$ 
given by Hara and Slade \cite{haraslade}, we also have the following corollary.
\begin{corollary}\label{corzd}
The critical value of the Bak-Sneppen model on $\mathbb{Z}^d$ satisfies
\[\frac{1}{2d+1}\leq p_c^{BS}(\mathbb{Z}^d)\leq \frac{1}{2d}+ 
\frac{1}{(2d)^2}+O\big(d^{-3}\big),
\qquad d\to\infty.\]
\end{corollary}

The paper is organised as follows. In Section \ref{sec2} we take some preliminary 
steps by describing an alternative way of constructing a Bak-Sneppen avalanche. 
Section \ref{sec3} uses this construction to  couple the Bak-Sneppen avalanche and 
another stochastic process. The proof that the critical value of the Bak-Sneppen 
avalanche is larger than that of the coupled stochastic process is given in Section 
\ref{sec4}. The proof of Theorem \ref{thm2} is completed in Section \ref{sec5} 
where we show that the coupled process  in fact constructs the cluster at the 
origin of site percolation with the origin always open.  In Section \ref{sec6}, we 
discuss some implications and generalisations of our methods and results.

\section{An alternative construction of the Bak-Sneppen model}\label{sec2}

In the introduction the Bak-Sneppen model was defined in its original format and 
then generalised to locally finite graphs. However, for our purposes it is more 
convenient to work with  an alternative construction of the Bak-Sneppen model. We 
call this new construction the {\em forgetful} Bak-Sneppen model, as the exact 
fitness values will be  no longer fixed  (or remembered). This idea borrows heavily 
from the `locking thresholds  representation' in \cite{zna1}, and was used in a 
much simpler form in \cite{gillett}.  The forgetful Bak-Sneppen model is defined 
below and then argued to be equivalent   to the normal Bak-Sneppen model,  in the 
sense that at all times,  the fitnesses have the same distributions.

Consider a Bak-Sneppen model on a finite transitive graph $G$ with $N$ vertices.  
To start with, all $N$ vertices have independent uniformly $(0,1)$ distributed 
fitnesses.  In the forgetful Bak-Sneppen model, all $N$ vertices have {\em fitness 
distributions}, instead of fitness values.  At time 0, all vertices have uniform 
(0,1) fitness distributions. The system at time $n$ is generated from the system at 
time $n-1$ by the following procedure.

\begin{enumerate}
\item We draw  $N$ new independent random variables according to the appropriate 
fitness distributions at time $n-1$.
\item The minimum of the fitnesses is found and fixed. 
\item All the other fitness values are discarded, and replaced by the 
conditional {\em distribution} of these fitnesses, given that they are larger than 
the 
observed minimal fitness. 
\item The vertex with minimal fitness and its neighbours have their value or 
fitness distributions replaced by  uniform (0,1) distributions. (So now all 
vertices have some distribution associated with them.)  This is the state of the 
system at time $n$.

\end{enumerate} 

It is easy to see that the fitness distributions at time $n$ generated  by this 
procedure are the same as the fitness distributions in the normal Bak-Sneppen model 
 at time $n$.

Furthermore,  all  fitness distributions have the convenient  property  that they 
are uniform distributions. Indeed, suppose that a  random variable $Y$ has a  
 uniform $(y,1)$ distribution, denoted by $F_y$. 
If we condition on $Y>z$, then $Y$ has 
distribution $F_{y\vee z}$, where $y\vee z=\max\{y,z\}$.  All our fitnesses 
initially  have uniform (0,1) distributions. Two things can change these  
distributions. They can be reset to $F_0$ by being updated, or they can be  
conditioned to be bigger than some given value; in both cases they remain uniform. 

The above construction gives a forgetful Bak-Sneppen  model on a finite graph, but 
this is easily extended to a forgetful Bak-Sneppen model on a locally finite graph. 
The only difference is that initially we assign the fitness distribution $F_0$  to 
the origin. The remaining vertices have a fitness distribution with all mass in the 
point 1,  denoted by $F_1$.  The avalanche ends when all the fitnesses within the 
avalanche are above the threshold, which is equivalent to saying that the minimal 
fitness is above the threshold. It is possible to see when the avalanche has 
finished by checking the value of the minimum fitness (phase 2 above). Thus we can 
use the forgetful method to generate avalanches.

\section{The construction of the coupling}\label{sec3}

This section is divided into three parts. To  begin with, some intuition behind the 
main result (Theorem \ref{thm2}) is given. This is followed by a  precise 
description of the coupling, and then  we give an example for added clarity.

\subsection{Intuition}

We are interested in comparing the Bak-Sneppen avalanche with the open cluster at 
the origin of independent site percolation, with the proviso that the origin
is open with probability 1 rather than with probability $p$. This clearly has no effect  on the critical value. 

Typically, site percolation is studied as a static random structure, but it is also 
possible to build up the open cluster at the origin dynamically. This is standard 
(we refer 
to \cite{grimmett} for details) but the idea can be described as follows. Starting with just
the origin, we can evaluate one of the neighbours and decide whether this neighbour
is open or not. If it is, we add it to the cluster, if it isn't, we declare it closed.
One can continue in this fashion, each time step evaluating neighbours of the current 
cluster one
by one. If the probability that a vertex is open, given the full history of
this process, is always equal to $p$, then in fact we do create the site-percolation open
cluster of
the origin. When there are no more unevaluated neighbours, the process stops,
and the cluster is finite in that case. 

The growth of both a Bak-Sneppen avalanche and the open cluster at  the origin is 
driven by the {\em extremal} vertices. In a Bak-Sneppen avalanche, the extremal 
vertices are those vertices that are contained within the avalanche and have 
neighbours outside the avalanche. It is only through one of the extremal vertices 
having the minimal fitness that the range of the avalanche can increase. For site 
percolation, the extremal vertices are those having a neighbour in the open cluster 
at the origin, but that are themselves unknown as to be open or closed. These are 
exactly the vertices at the edge of the cluster, and they  will increase the size 
of the cluster by being open. Since it is the extremal vertices that drive the 
spread of both processes, the task is to relate the two sets of extremal vertices 
to each other.

The major difficulty to overcome is that in the Bak-Sneppen model an extremal 
vertex may be updated  by neighbouring activity before having  minimal fitness 
itself, whereas in site percolation a vertex is either open or closed. So in the  
Bak-Sneppen model it is possible that a previously active extremal vertex never has 
minimal fitness, having been made inactive by a subsequent neighbouring update. 
Conversely,  an originally inactive vertex can be  made active. Hence, in the 
Bak-Sneppen model the neighbour of an active vertex will not necessarily be 
updated, while in our construction of the open cluster at the origin in site 
percolation, the neighbour of an open site is always considered.  This means that 
it is not useful to couple the two models in the natural manner by realising the 
fitness and determining if the vertex is open and closed immediately with the same 
random variable.

The following heuristics make Theorem \ref{thm2} plausible. If  a vertex's fitness 
is not minimal, then its conditional 
distribution based on this information is stochastically larger than its original 
uniform $(0,1)$ distribution. So if a vertex is updated by a neighbour having 
minimal fitness, this makes its fitness stochastically smaller, making the vertex  
more likely to be active and therefore, intuitively at least, the avalanche is  
more likely to continue. This means that on average the interference from the 
non-extremal vertices of the Bak-Sneppen model on the extremal  vertices should be 
beneficial to the spread of the avalanche.

\subsection{The coupling}\label{sec:coupling}

We now describe the construction of a process that we will refer to as the {\em 
coupled process}.  As we shall see later, this process is constructed  
 in such a way that it is stochastically dominated by 
the Bak-Sneppen avalanche, which  is crucial for our argument.
In Section \ref{sec5} we show that this coupled process in 
fact  constructs  the cluster at the origin of site percolation. 

Let $V(G)$ be the vertex set of the graph $G$. The coupled process is a stochastic 
process with values in $\{([0,1]\times \{f, d\})^A, A\subset V(G)\}$. An entry $(a, 
f)$ means that the value of that vertex is fixed at $a$ forever, while an entry 
$(a, d)$ means that the value of that vertex is distributed uniformly on $(a,1)$. 
The coupled process is coupled to a forgetful  Bak-Sneppen avalanche, and is 
constructed as follows.

Fix an avalanche threshold $p$. We start with two copies of the graph $G$, denoted 
by $G_B$ (for the Bak-Sneppen avalanche) and $G_C$ (for the coupled process). 
Initially we assign the value $0$ to the origin of $G_B$ and $(0,f)$ to the origin 
of $G_C$, and we call the origin in $G_C$ {\em open} (as anticipated before). 
Then all the $\Delta$ neighbours of the 
origin of both graphs get distribution $F_0$. On $G_C$,  we define the {\em 
extremal set} ${\cal E}$ as the set of all points that have been assigned a 
distribution, but not (yet) an exact value. 

The Bak-Sneppen avalanche on $G_B$ is generated according to the aforementioned 
(forgetful) construction, i.e., we sample new fitnesses,  fix the minimal value and 
then calculate the fitness distributions accordingly.  In the coupled process, only 
the vertices contained in ${\cal E}$ are  considered. We apply the following 
procedure to all vertices in ${\cal E}$.

Consider a vertex $v_C\in {\cal E}$  with  $G_B$-counterpart $v_B$. Let $F_z$ and 
$F_y$ be their respective fitness  distributions. We realise the fitnesses of the 
vertices in $G_B$, and in particular realise the fitness of $v_B$ with an 
independent uniformly (0,1) distributed random variable $U$ via $y + (1-y)U$. Let 
$M$ be the minimal fitness in $G_B$. As long as the vertex with minimal fitness in 
the avalanche is active, i.e., $M \leq p$,  we have the following two options, with 
corresponding rules  for the coupled process.  One should bear in mind that 
 the main goal of the coupling is the stochastic domination. In Section \ref{sec:example} 
below, these are illustrated by  an explicit example.

\begin{enumerate}

\item The fitness of $v_B$ is not minimal.\\
We alter the distribution of $v_C$ by conditioning on the extra information that 
the fitness of $v_B$ must be bigger than $M$. Since the fitness of $v_B$ is not 
minimal, we have $y+(1-y)U>M$, and hence $U>(M-y)^+/(1-y)$.  The new distribution 
of $v_C$ is $F_{\hat{z}}$, where
\[\hat{z}=z+(1-z)\frac{(M-y)^+}{(1-y)}.\]

\item The fitness of $v_B$ is minimal, so it has value $M$.\\
It follows that $y+(1-y)U=M$. The fitness of $v_C$ is now {\em fixed} at 
\[z+(1-z)U=z+(1-z)\frac{(M-y)}{(1-y)}.\]
If this value is less than $p$, we say that $v$ is {\em open}, remove $v$ from  
${\cal E}$, add the neighbours of $v$ that have an undetermined state to ${\cal  
E}$, and give them distribution $F_0$. If the value of $v$ is larger than $p$,  
then $v$ is {\em closed} and removed from ${\cal E}$.

\end{enumerate}

The final step of the construction is as follows.  The first time that the vertex with 
minimal fitness in $G_B$ is inactive (that is,  $M>p$), the Bak-Sneppen avalanche 
has finished. As soon as this happens, we fix {\em all}  the values of the vertices 
in ${\cal E}$  in the following way,  similar  to rule $2$ above.  Let  $v_C \in 
{\cal E}$ and  $v_B$ have  fitness distributions $F_z$ and $F_y$ respectively, and 
let $U$ be  the associated uniform (0,1) random variable. Then $U$ satisfies  
$y+(1-y)U\geq M$, i.e., $U\geq  (M-y)/(1-y)$. The new distribution of  $v_C$ is  
$F_{\hat{z}}$ with $\hat{z}=z+(1-z)(M-y)^+/(1-y)$. As final step of the coupling, 
we realise the fitness of $v_C$ as $\hat{z}+(1-\hat{z})X$, where $X$ is an 
independent uniformly (0,1) distributed random variable. In Section \ref{sec4} we 
show that as soon as the Bak-Sneppen avalanche ends, this fitness value is at least 
$p$, and hence all the vertices in ${\cal E}$ will be  closed. Before that, we give 
an example to illustrate the coupling procedure described above.

\subsection{An example}\label{sec:example}

\begin{figure}
\begin{picture}(400,400)(-50,0)
\put(-30, 370){$a)$}
\put(200, 370){$b)$}
\put(-30, 270){$c)$}
\put(200, 270){$d)$}
\put(100, 360){\line(0,1){20}}
\put(300, 360){\line(0,1){20}}
\put(100, 360){\line(0,1){20}}
\put(300, 360){\line(0,1){20}}
\put(70, 340){\circle*{4}}
\put(130, 340){\circle*{4}}
\put(100, 360){\line(-3,-2){30}}
\put(100, 360){\line(3,-2){30}}
\put(100, 360){\circle*{4}}
\put(135, 340){$F_{0}$}
\put(45, 340){$F_{0}$}
\put(300, 360){\circle*{4}}
\put(107, 360){$F_0$}
\put(270, 340){\circle*{4}}
\put(330, 340){\circle*{4}}
\put(300, 360){\line(-3,-2){30}}
\put(300, 360){\line(3,-2){30}}
\put(300, 360){\circle*{4}}
\put(245, 340){$U_2$}
\put(335, 340){$U_3$}
\put(307, 360){$U_1$}
\put(70, 240){\circle*{4}}
\put(70, 240){\circle{7}}
\put(130, 240){\circle*{4}}
\put(100, 260){\line(-3,-2){30}}
\put(100, 260){\line(3,-2){30}}
\put(100, 260){\circle*{4}}
\put(100, 260){\line(0,1){20}}
\put(135, 240){$F_M$}
\put(50, 240){$M$}
\put(107, 260){$F_M$}
\put(300, 260){\line(0,1){20}}
\put(300, 260){\circle*{4}}
\put(270, 240){\circle*{4}}
\put(330, 240){\circle*{4}}
\put(250, 220){\circle*{4}}
\put(290, 220){\circle*{4}}
\put(300, 260){\line(-3,-2){30}}
\put(300, 260){\line(3,-2){30}}
\put(270, 240){\line(-1,-1){20}}
\put(270, 240){\line(1,-1){20}}
\put(300, 260){\circle*{4}}
\put(335, 240){$F_M$}
\put(245, 240){$F_0$}
\put(230, 220){$F_0$}
\put(295, 220){$F_0$}
\put(307, 260){$F_0$}
\put(-30, 150){$e)$}
\put(200, 150){$f)$}
\put(-30, 50){$g)$}
\put(200, 50){$h)$}
\put(100, 160){\line(0,1){19}}
\put(300, 160){\line(0,1){19}}
\put(100, 60){\line(0,1){19}}
\put(300, 60){\line(0,1){19}}
\put(70, 140){\circle{4}}
\put(130, 140){\circle{4}}
\put(100, 160){\line(-3,-2){29}}
\put(100, 160){\line(3,-2){29}}
\put(100, 160){\circle*{4}}
\put(135, 140){$F_{x}$}
\put(50, 140){$F_{x}$}
\put(300, 160){\circle*{4}}
\put(100, 60){\circle*{4}}
\put(270, 140){\circle{4}}
\put(330, 140){\circle{4}}
\put(300, 160){\line(-3,-2){29}}
\put(300, 160){\line(3,-2){29}}
\put(300, 160){\circle*{4}}
\put(335, 140){$F_x$}
\put(335, 125){{\tt rule 1}}
\put(250, 140){$F_{x}$}
\put(250, 125){{\tt rule 2}}
\put(130, 40){\circle{4}}
\put(100, 60){\line(3,-2){29}}
\put(135, 40){$F_{\hat{x}}$}
\put(300, 60){\circle*{4}}
\put(270, 40){\circle*{4}}
\put(330, 40){\circle{4}}
\put(250, 20){\circle{4}}
\put(290, 20){\circle{4}}
\put(300, 60){\line(-3,-2){29}}
\put(300, 60){\line(3,-2){29}}
\put(270, 40){\line(-1,-1){19}}
\put(270, 40){\line(1,-1){19}}
\put(300, 60){\circle*{4}}
\put(335, 40){$F_{\hat{x}}$}
\put(295, 20){$F_0$}
\put(230, 20){$F_0$}
\put(270, 140){\circle{7}}
\end{picture}
\caption{
Graphs $a)$ -- $d)$ show a time step in the forgetful  Bak-Sneppen process, and  
graphs $e)$ -- $h)$ show a time step in the coupled process. The encircled vertex 
has the minimal fitness. In the coupled process, black points are open, white  
points are undetermined, and closed points are omitted. }
\label{fig3}
\end{figure}
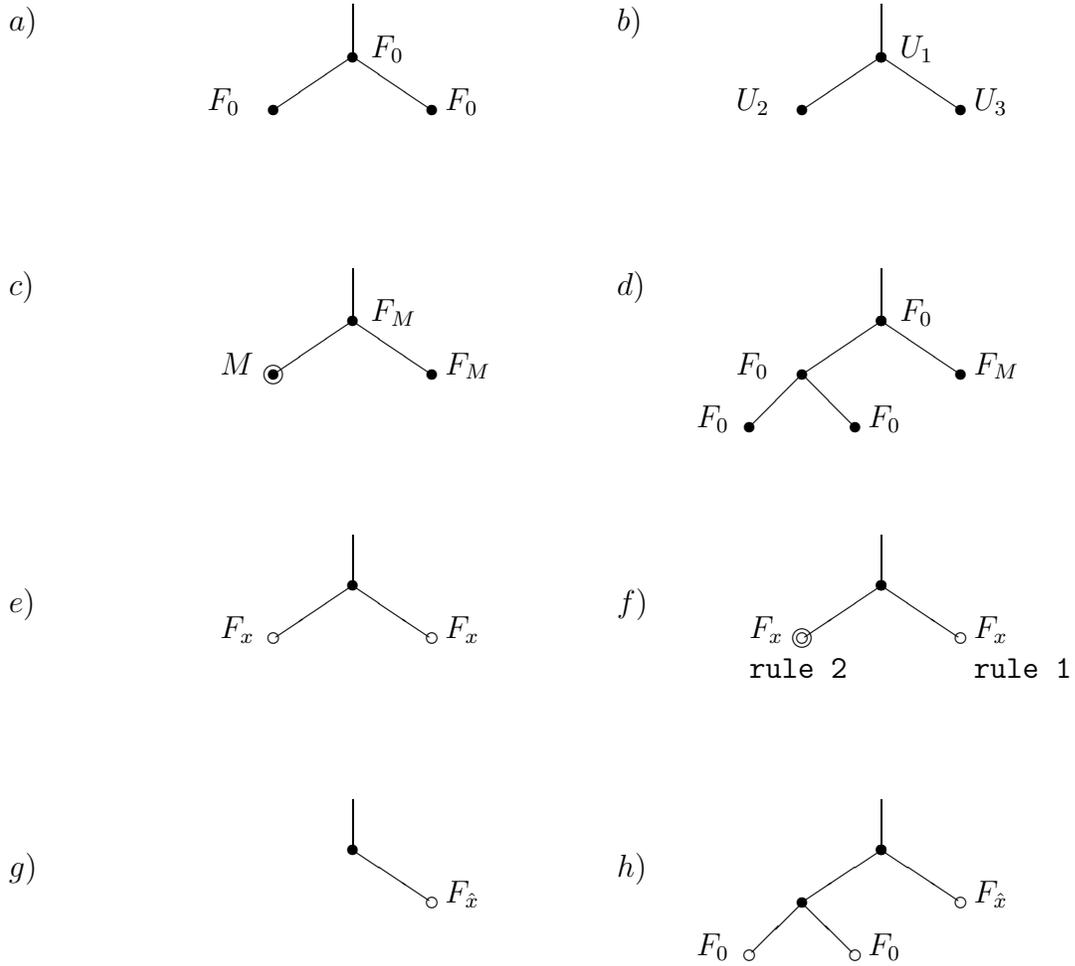

The behaviour of the processes is illustrated by the following example, displayed  
in Figure \ref{fig3}. In this example the graph $G$ is a tree. For illustration 
purposes, we show only the part of the graph where the activity takes place.

Consider the forgetful avalanche at time $n$ say, in the following situation, see 
Figure \ref{fig3}, graph $a$: all fitnesses shown have distribution $F_0.$  Then  
in graph $b$, uniform $(0,1)$ random variables are drawn. The random variables 
$U_1, U_2,$ and $ U_3 $ are associated with the vertices visible in the picture. 
Other random variables are of course drawn for the other vertices. We then use the  
full set of random variables to determine the location and the magnitude of the new 
minimal fitness. This happens to be the vertex corresponding to the random variable 
$U_2$ (graph $c$). Finally, the  new fitness distributions for time $n+1$ are 
determined (graph $d$). Note that the forgetful Bak-Sneppen  model never actually 
assumes the values given in graph $b$. 

During the same time step, the coupled process evolves as follows. We only consider 
the vertices in the extremal set ${\cal E}$, see Figure \ref{fig3}, graph $e$. 
Before the time step, the vertices have fitness distributions $F_x$, for some $x\in 
(0,1)$. Given the location of the minimal fitness in the Bak-Sneppen avalanche, the 
vertices in  ${\cal E}$ are classified according to the rules 1 and 2 above (graph 
$f$).  From the location and magnitude of the minimal fitness of the avalanche, it 
follows that  $U_3 \geq M$, so $\hat{x}= x+(1-x)M$. Finally, the value of the 
vertex that corresponds with the vertex with minimal fitness in the avalanche is 
fixed, according to rule 2.  Its value $f$ is given by  $f=x+(1-x)M$. Now there are 
 two possible cases: either $f>p$, and the vertex is closed (graph $g$),  or $f\leq 
p$, and the vertex is open and its undetermined neighbours are added to ${\cal E}$ 
with  distribution  $F_0$ (graph $h$). 

\section{A domination principle}\label{sec4}

To show that the critical value of the coupled process can be no smaller  than the 
critical value of the Bak-Sneppen avalanche, we use a domination argument. The 
propositions below show that the coupled process can finish no later than the 
Bak-Sneppen avalanche (so that the avalanche can be said to dominate the coupled 
process). 
\begin{proposition}\label{prop5}
For every $v_C \in G_C$ and corresponding $v_B \in G_B$, at all times, the 
(conditional) fitness distribution of $v_C$ is stochastically larger than the 
(conditional) fitness distribution of $v_B$. 
\end{proposition}
\noindent {\bf Proof:}
It should be noted that this proposition only makes sense for vertices in ${\cal 
E}$. Furthermore, it is safe to assume that the $p$-avalanche is still in progress, 
so the minimal fitness is less than $p$. The proof proceeds by induction. When new 
vertices are added to the coupled process, they (and their equivalents in $G_B$) 
have uniform $(0,1)$  distributed fitnesses. This is by definition for the coupled 
process, but also  holds for $G_B$, since  vertices in $G_B$ corresponding to new 
vertices added to ${\cal E}$ are always neighbours of the vertex with minimal 
fitness. This means that the statement of the proposition holds for new vertices 
added to ${\cal E}$.

To make the induction step,  consider $v_C \in {\cal E}$ with corresponding vertex 
$v_B$ in $G_B$,  and let $F_z$ and $F_y$ be  the fitness distributions of $v_C$ and 
$v_B$ at time $n$, where $y\leq z<1$.  Let $u$ be the realisation of the uniform 
$(0,1)$ random variable associated with $v_C$ and $v_B$ at the intermediate step, 
and let $m$ be the minimal  fitness. 

Assume first that $v_B$ does not have the minimal fitness. This provides 
information on the value of $u$, namely that $y+(1-y)u > m$, and hence 
$u>(m-y)/(1-y)$. If $y\geq m$, this information is useless: we already knew that 
$u>0$, and the fitness distributions of $v_C$ and $v_B$ are not changed. If $y < 
m$, then the inequality for $u$ does contain information, and we can calculate the 
corresponding inequality for the fitness of $v_C$:
\[ z+(1-z)u  >  z+\frac{(1-z)(m-y)}{1-y}  =  m + \frac{(1-m)(z-y)}{1-y}:=\hat{z}.\]
So at time $n+1$, $v_B$ has distribution $F_{y \vee m}$ and $v_C$ has  distribution 
$F_{\hat{z}}$. Since $m,y<1$ and $y\leq z$, we have $\hat{z}\geq m$. Hence $(y \vee 
m) \le  \hat{z}$,  and the desired property holds. 

Second, we consider the case that a neighbour of $v_B$ had minimal fitness. In that 
case the fitness distribution of $v_B$ is reset to $F_0$, and there is nothing left 
to prove. $\hfill \Box$

\begin{proposition}\label{propE}
At the moment that the $p$-avalanche ends, all vertices in ${\cal E}$ are closed.  
As a consequence, if the probability of an infinite $p$-avalanche is zero, then 
there cannot be an infinite cluster of open sites in the coupled process, almost 
surely.
\end{proposition}
{\bf Proof:}
By Proposition \ref{prop5}, at all times every point in  ${\cal E}$ has a fitness 
that is stochastically larger than the fitness of the corresponding vertex in the 
avalanche. Hence, if the $p$-avalanche ends, then in the coupled process all 
neighbours in the set ${\cal E}$ will be closed, as their fixed values can not be 
smaller than those in the avalanche, which are already greater than $p$ as the 
avalanche has ended. This removes all vertices from ${\cal E}$ and ensures that no 
more are added, implying that in  the coupled process  no more vertices will be 
added to the open cluster around the origin. $\hfill \Box$\\ 

 We conclude this section by giving an example where the coupled process is finite, but the Bak-Sneppen avalanche is infinite. This shows that the stochastic domination described in this section is not a stochastic equality.

Let $G=\mathbb{Z}$ and  $p=0.7$. Suppose that both in the first step and the second step in the Bak-Sneppen model, the origin is minimal with fitness $0.5$. In the coupled process, the neighbours of the origin have  fitness distribution  $F_{0.5}$ after the first step,
and $F_{0.5+(1-0.5)0.5}=F_{0.75}$ after the second step. Since $0.75>p$, this implies that the neighbours of the origin will eventually be closed, and the cluster in the coupled process is finite. However, the Bak-Sneppen avalanche may very well be infinite. 

\section{The  cluster at the origin of site percolation}\label{sec5}

To complete the proof of Theorem \ref{thm2}, it remains to show that the coupled 
process  in fact constructs the open cluster at the origin of independent site 
percolation, with the proviso that the origin is open with probability 1. To get into 
the right frame of mind for 
the proof, we first give an example. At the same time, the example illustrates the 
construction of the coupled process in action.

\subsection{An example}  Consider the Bak-Sneppen avalanche and the coupled process 
defined on $\mathbb{Z}$ with parameter $p$. We wish to calculate the probability 
that in the coupled process both neighbours of  the origin are closed. Note that 
for the site percolation cluster this probability is $(1-p)^2$, so our aim is to 
show that this probability is also $(1-p)^2$ for the coupled process. 
To calculate this probability, we introduce the following, more general probability: 
for all  $0\leq x\leq p$, let $g_p(x)$ be the probability that both neighbours of the 
origin will be declared closed, given that their current fitness 
distributions both are $F_x$. In this notation, the desired probability is equal to 
$g_p(0)$.

Starting with the distributions $F_x$ for both neighbours, we call the
first subsequent step, the {\em first} time step.
For the coupled process, both neighbours of the origin are declared closed if their 
realised 
values are above $p$. Noting that both neighbours have distribution $F_x$, this will 
happen at the first time step if in the 
Bak-Sneppen model all three values are above $(p-x)/(1-x)$. 
If the minimum, which has density $3(1-b)^2$, is below $(p-x)/(1-x)$, and located 
at the origin (which happens with probability 1/3), then we have to look at 
subsequent updates in the Bak-Sneppen model.  

In this second case,  the three 
fitness distributions in the Bak-Sneppen model are reset to $F_0$. However, in the 
coupled process, the fitnesses of $-1$ and $1$ are  now $F_{x+(1-x)b}$, where $b$ 
is the avalanche minimum at the first time step. 
For the second time step, we are now in a similar situation as for the first,  
except that the fitness distribution has a different parameter: $x+(1-x)b$ instead 
of $x$. This similarity holds for any starting level $x$,
and leads to the following expression for $g_p(x)$:
\begin{align}\label{marku}
g(x): =g_p(x)&=\Big(\frac{1-p}{1-x}\Big)^3+\frac{1}{3}\int_0^{\frac{p-x}{1-x}}
3(1-b)^2 g(x+(1-x)b)d b.\end{align}
Substituting $y=x+(1-x)b$, equation (\ref{marku})  becomes
\begin{equation}\label{gx} g(x) =\frac{1}{(1-x)^3}\Big((1-p)^3+ \int_x^p(1-y)^2 
g(y)d y   \Big).\end{equation}
Using (\ref{gx}), a little algebra yields that for   small $h$,
\begin{equation}\label{gxh} g(x+h)-g(x)
=\frac{(1-x)^3-(1-x-h)^3}{(1-x-h)^3}g(x)-\frac{1}{(1-x-h)^3}\int_x^{x+h}(1-y)^2 
g(y)d y.\end{equation}
 Since $0\leq g\leq 1$, it follows from (\ref{gxh}) that $g(x+h)-g(x)\to 0$ for $h\to 0$, so $g$ is  continuous.
Hence, we can calculate the differential quotient:  
\[\lim_{d\downarrow  0}  \frac{g(x+h)-g(x)}{h}=\frac{3(1-x)^2g(x)}{(1-x)^3} -\frac{(1-x)^2g(x)}{(1-x)^3}=
\frac{2g(x)}{1-x}.\]
 The same holds for the left-hand limit, so 
$g(x)$ is differentiable, and $g'(x)=2g(x)/(1-x)$. This differential 
equation has a unique solution for each $p$, given by
$g(x)=c(p)/(1-x)^2$.
Using the boundary condition $g_p(p)=1$, we find $c(p)=(1-p)^2$, so that
\[ g_p(x)=\frac{(1-p)^2}{(1-x)^2}.\]
In particular, the  desired probability that in the coupled process both neighbours 
are closed is given by $g_p(0)=(1-p)^2$, as required.

Although this example gave us what we wanted,  clearly  this type of calculation 
does not generalise to  more complicated events. Therefore, the proof that the 
coupled process constructs the site percolation open cluster, which we turn to now, 
necessarily has a different flavour.

\subsection{The proof}

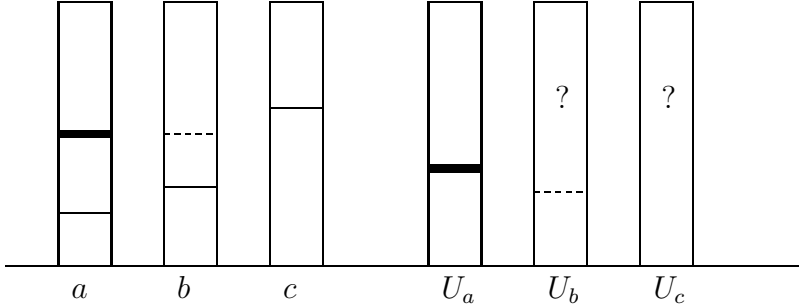
\begin{figure}[h!]
\begin{center}
\begin{picture}(300,130)(0,-20)

\put(-20, 0){\line(1,0){300}}
\put(0, 0){\line(0,1){100}}
\put(20, 0){\line(0,1){100}}
\put(0, 100){\line(1,0){20}}
\put(0, 20){\line(1,0){20}}
\put(0, 50){\linethickness{1mm}\line(1,0){20}}
\put(5, -12){$a$}

\put(40, 0){\line(0,1){100}}
\put(60, 0){\line(0,1){100}}
\put(40, 100){\line(1,0){20}}
\put(40, 30){\line(1,0){20}}
\put(40, 50){\line(1,0){2}}
\put(44, 50){\line(1,0){2}}
\put(48, 50){\line(1,0){2}}
\put(52, 50){\line(1,0){2}}
\put(56, 50){\line(1,0){2}}
\put(45, -12){$b$}

\put(80, 0){\line(0,1){100}}
\put(100, 0){\line(0,1){100}}
\put(80, 100){\line(1,0){20}}
\put(80, 60){\line(1,0){20}}
\put(85, -12){$c$}

\put(140, 0){\line(0,1){100}}
\put(160, 0){\line(0,1){100}}
\put(140, 100){\line(1,0){20}}
\put(140, 37){\linethickness{1mm}\line(1,0){20}}
\put(145, -12){$U_a$}

\put(180, 0){\line(0,1){100}}
\put(200, 0){\line(0,1){100}}
\put(180, 100){\line(1,0){20}}
\put(180, 28){\line(1,0){2}}
\put(184, 28){\line(1,0){2}}
\put(188, 28){\line(1,0){2}}
\put(192, 28){\line(1,0){2}}
\put(196, 28){\line(1,0){2}}
\put(185, -12){$U_b$}
\put(188, 60){$?$}

\put(220, 0){\line(0,1){100}}
\put(240, 0){\line(0,1){100}}
\put(220, 100){\line(1,0){20}}
\put(220, 0){\line(1,0){20}}
\put(225, -12){$U_c$}
\put(228, 60){$?$}

\end{picture}
\caption{In the forgetful BS-avalanche, before the update, vertices $a$, $b$ and 
$c$ have fitness distributions  $F_{0.2}$,  $F_{0.3}$, and  $F_{0.6}$, 
respectively.  After realising these distributions, vertex $a$ is minimal with 
value $M=0.5$.  This means that $U_a=\frac{0.5-0.2}{1-0.2}=\frac{3}{8}$,  $U_b\geq 
\frac{0.5-0.3}{1-0.3}=\frac{2}{7}$, and $U_c\geq 0$. This sample, namely $U_a=3/8$ 
combined with a sample from a uniform $(2/7,1)$ and a uniform $(0,1)$ distribution, 
is a sample of three i.i.d.\ uniform $(0,1)$ random variables.}
\label{fig4}
\end{center}
\end{figure}

Our first goal is to determine the distribution of the information we use to 
generate the coupled process. More precisely, consider an arbitrary step of the 
forgetful Bak-Sneppen model, when there are $n$ vertices in the avalanche range so 
far. We enumerate these vertices $1,\ldots, n$, and suppose that all $n$ vertices 
in the avalanche have just been assigned a (conditional)  distribution 
$F_{y_1},\ldots, F_{y_n}$. (Recall that these are just uniform distributions above 
the respective $y_i$'s.) We sample from this random vector, using independent 
uniform $(0,1)$ distributed random variables  $U_1,\ldots, U_n$: a sample from 
$F_{y_i}$ is realised via $y_i + (1-y_i)U_i$. We locate the minimum $M$, at vertex 
$K$ say; note that both $M$ and $K$ are random. Hence,
\begin{equation}\label{verg}
U_K =\frac{M-y_K}{1-y_K}.
\end{equation}
Conditional on $K$ and $M$, the remaining values $U_i$,  $i \neq K$, are  uniformly 
distributed above $\max\{y_i, M\}$ respectively, that is, we know that
$$
U_i > \frac{(M-y_i)}{1-y_i}, \qquad i \neq K.
$$ 
When we now also sample from all the other entries $i \neq K$, (which are uniformly 
distributed above $(M-y_i)^+/(1-y_i)$ respectively) we have described a somewhat 
complicated 
way of sampling from the original vector $(U_1,\ldots, U_n)$, that is, such a 
sample yields independent uniform $(0,1)$ distributed entries, see also Figure 
\ref{fig4} and its caption. Note that we do not claim that $U_K$ is uniformly 
distributed on $(0,1)$: it is not. However, since the index $K$ is random, 
this does not 
contradict the fact that the vector $(U_1,\ldots, U_n)$ consists of independent uniform 
$(0,1)$ random variables.

Looking back to Section 3.2, it should be clear that in the coupled process  
independent uniform $(0,1)$ random variables generated in the above way are used to 
alter the values of the vertices  contained in ${\cal E}$. Note that  using $|{\cal 
E}|$ entries rather than $n$  does not affect their marginal  distributions or 
dependence structure, as the values of $U_i's$ do not depend on whether the 
associated vertices are in ${\cal E}$ or not. 

It is now possible to give a direct description of the construction of the coupled 
process. We start with the origin being open and look at the neighbours of the 
origin, which initially have distribution $F_0$. These distributions are realised 
as follows, using the independent uniform $(0,1)$ random variables described above. 
At each time step at most one value becomes fixed and the rest are given 
distributions. The fixed value corresponds to the case that $K \in {\cal E}$. To  
calculate the new values of vertices in ${\cal E} \backslash \{K\}$, we use the 
information 
that the $U_i$'s are independently and uniformly distributed above 
$(M-y_i)^+/(1-y_i)$. This means that we do not fix their actual values at that time 
step, but instead change their distributions conditioned on this information. Once 
a vertex has a fixed value, it is declared open if and only if this value falls 
below $p$. Whenever a vertex is declared open,  the neighbours that neither have a 
fixed value nor belong to ${\cal E}$ are added to ${\cal E}$  with  distribution 
$F_0$. 

Since fitnesses are initially independent uniform $(0,1)$ when 
added to ${\cal E}$ and the information we use to update the distributions is also 
independent uniform $(0,1)$, the following holds:  if at any time point the 
procedure is stopped and all the distributions are realised, one will recover an 
independent uniform $(0,1)$ sample.  Hence, all considered vertices (except the 
origin) are open independently  and with probability $p$. It should now be obvious 
that our procedure is  no different to building a site percolation cluster at the 
origin by the iterative method of assigning independent uniform $(0,1)$ random 
variables  to all undetermined neighbours of the cluster and calling a vertex open 
if  its random variable takes a value less than $p$.  This completes the proof of  
Theorem \ref{thm2}.~$\hfill \Box$\\

\noindent
Note that in case of an  infinite Bak-Sneppen avalanche,  some vertices in the 
coupled process may never get a fixed value. This is not a problem, because this is 
just what happens if an infinite open cluster around the origin is built up 
dynamically: not all vertices will be tested in the process of constructing this
cluster.

\section{Final remarks and extensions}\label{sec6}
Throughout this paper we have only considered locally finite transitive graphs.  We 
assumed transitivity to avoid technicalities that would have obscured the main  
lines of reasoning. However,  our results also hold in a more general  setting, 
namely for any locally finite graph. The following observations  explain this 
generalisation.
The lower bound (Proposition \ref{thm1}) can easily be adapted  by considering a 
branching process with binomial$(\Delta^*+1,p)$  offspring,  where $\Delta^*$ is 
the maximal degree of the graph. Note that the lower bound is trivial if  
$\Delta^*=\infty$. 
 
The coupling argument  used to prove that the Bak-Sneppen avalanche dominates site  
percolation, at no point used the transitivity of the underlying graph, and hence  
also holds for non-transitive graphs. However, for non-transitive graphs, the 
choice  of the  origin affects the behaviour of the avalanche. The upper bound 
(Theorem \ref{thm2}) is generalised by the following observation: although the 
distribution of the size of the open cluster around the origin in site percolation 
does depend on the choice of the origin, standard arguments yield that the critical 
value does not.

Another consequence of our methods is the following.
The careful reader may have noticed that the proofs actually yield a stronger 
result than stated in Theorem \ref{thm2}, namely {\em stochastic domination}. 
Define the range of site percolation to be the cardinality of the open cluster 
around the origin plus all its closed neighbours (these closed neighbours 
correspond to updated vertices in the Bak-Sneppen avalanche that were never 
minimal). The proof of Theorem \ref{thm2} then demonstrates that the range of the 
$p$-avalanche  is {\em stochastically larger} than the range of site percolation 
with parameter $p$. 

Although not explicitly stated in the proof of Proposition \ref{thm1}, a similar 
extension also applies there. The set of offspring of a branching process with a 
binomial $(n-1,p)$ offspring distribution is equivalent to the open cluster around 
the origin (root) of site percolation with parameter $p$ on 
$T_n^*$, where $T_n^*$ is a rooted tree  where the root has degree $n-1$, and all 
other vertices have degree $n$. In this case we get that the range of a 
$p$-avalanche on a transitive graph with common vertex degree $\Delta$ is {\em 
stochastically smaller} than the range of site percolation on $T_{\Delta+2}^*$.

Finally, we argue that Theorem \ref{thm2} holds as well for the critical value 
(\ref{crit2}). It is well-known that for site percolation on  $\mathbb{Z}^d$ or on 
a tree, $p_c^{site}(G)$ is equal to the critical value associated with the expected 
size of the open cluster at the origin, see Grimmett \cite{grimmett}. Since each 
vertex in the open cluster  contributes at most $\Delta$ closed neighbours to the 
range of site percolation,   the range  is always less than $\Delta$ times the size 
of the cluster. Hence, the critical values associated with the expectation of these 
two objects are the same. As a consequence, the stochastic bounds given above imply 
 that the bounds in Proposition \ref{thm1} and Theorem \ref{thm2} also hold for the 
critical value (\ref{crit2}).\\ \\ \\
\noindent
{\bf Acknowledgements}\\
The authors would like to thank the referee for his or her comments, which have 
improved the presentation and readability of the paper.

\bibliographystyle{plain}

\end{document}